\documentclass[10pt]{amsart}
\usepackage{color}
\usepackage{verbatim}
\usepackage[pdftex]{graphicx}

\usepackage{amssymb}
\allowdisplaybreaks[4]

\def\LT{\left}
\def\RT{\right}

\definecolor{c20}{rgb}{0.,0.7,0.}
\definecolor{c30}{rgb}{0.,0.,1.}
\definecolor{c40}{rgb}{1,0.1,0.7}
\definecolor{c50}{rgb}{1,0,0}
\definecolor{c60}{rgb}{1,0.9,0.1}

\newcommand{\abs}[1]{\left\lvert #1 \right\rvert}

\newcommand{\pk}[1]{\mathbb{P} \left\{ #1 \right \} }

\newcommand{\R}{\mathbb{R}}

\newcommand{\BQN}{\begin{eqnarray}}
\newcommand{\EQN}{\end{eqnarray}}
\newcommand{\BQNY}{\begin{eqnarray*}}
\newcommand{\EQNY}{\end{eqnarray*}}

\newcommand{\BS}{\begin{sat}}
\newcommand{\ES}{\end{sat}}
\newcommand{\BT}{\begin{theo}}
\newcommand{\ET}{\end{theo}}
\newcommand{\BK}{\begin{korr}}
\newcommand{\EK}{\end{korr}}

\newcommand{\BD}{\begin{de}}
\newcommand{\ED}{\end{de}}

\newcommand{\BIT}{\begin{itemize}}
\newcommand{\EIT}{\end{itemize}}
\newcommand{\BDI}{\begin{description}}
\newcommand{\EDI}{\end{description}}

\newcommand{\BRM}{\begin{remarks}}
\newcommand{\ERM}{\end{remarks}}

\newcommand{\BEL}{\begin{lem}}
\newcommand{\EEL}{\end{lem}}

\newtheorem{theo}{Theorem}[section]
\newtheorem{sat}[theo]{Proposition}
\newtheorem{de}[theo]{Definition}
\newtheorem{lem}[theo]{Lemma}

\newtheorem{korr}[theo]{Corollary}
\newtheorem{remark}[theo]{Remark}
\newtheorem{remarks}[theo]{Remarks}
\newtheorem{prop}[theo]{Proposition}

\newcommand{\nelem}[1]{{Lemma \ref{#1}}}

\newcommand{\COM}[1]{}

\newcommand{\QED}{\hfill $\Box$}

\newcommand{\kb}[1]{\boldsymbol{#1}}
\newcommand{\vk}[1]{\kb{#1}}

\def\IF{\infty}

\def\LT{\left}
\def\RT{\right}

\def\K1#1{\textcolor{cyan}{#1}}

\def\ccj#1{\textcolor{black}{#1}}

\def\K1#1{\textcolor{cyan}{#1}}

\def\jc#1{\textcolor{black}{#1}}
\def\kk#1{\textcolor{black}{#1}}

\def\Lcc#1{\textcolor{black}{#1}}

\def\oA{\overline A}
\def\oB{\overline B}

\begin{document}

\title{Logarithmic asymptotics for probability of component-wise ruin in a two-dimensional Brownian model}
\author{Krzysztof D\c{e}bicki}
\address{Krzysztof D\c{e}bicki, Mathematical Institute, University of Wroc\l aw, pl. Grunwaldzki 2/4, 50-384 Wroc\l aw, Poland}
\email{Krzysztof.Debicki@math.uni.wroc.pl}

\author{Lanpeng Ji}
\address{Lanpeng Ji, School of Mathematics, University of Leeds, Woodhouse Lane, Leeds LS2 9JT, United Kingdom
}
\email{l.ji@leeds.ac.uk}

\author{Tomasz Rolski}
\address{Tomasz Rolski, Mathematical Institute, University of Wroc\l aw, pl. Grunwaldzki 2/4, 50-384 Wroc\l aw, Poland}
\email{Tomasz.Rolski@math.uni.wroc.pl}


\bigskip

\date{\today}
 \maketitle

 {\bf Abstract:}

We consider a two-dimensional ruin problem where the surplus process
of business lines
is modelled by a two-dimensional correlated Brownian motion with drift.
We study the ruin function $P(u)$ for the component-wise ruin (that is both  business lines are ruined in an infinite-time horizon),
where $u$ is the same initial capital for each line.
We measure the goodness of the business by analysing  the adjustment coefficient, that is the
limit of $-\ln P(u)/u$ as $u$ tends to infinity,
which depends essentially on the correlation $\rho$ of the two surplus processes.
In order to work out the adjustment coefficient we
solve a two-layer optimization problem.

 {\bf Key Words:}  Adjustment coefficient;  logarithmic asymptotics; quadratic programming problem; ruin probability; two-dimensional Brownian motion
 
 {\bf AMS Classification:} Primary 60G15; secondary 60G70

 \section{Introduction}

 In classical risk theory, the surplus process of an insurance company is modelled by the compound Poisson risk model.
 For both applied and theoretical investigations, calculation of ruin probabilities for such model is of particular interest.
 In order to avoid technical calculations,  {\it diffusion approximation} is often considered
 (e.g., \cite{Igl69, Gran91, KPW12, AsmAlb10}),
 which results in \kk{tractable} approximations for the interested finite-time or infinite-time ruin probabilities. The basic premise for the approximation is to let the number of claims grow in a unit time interval and to make the claim sizes smaller in such a way that the risk process converges to a Brownian motion with drift. Precisely, the Brownian motion  risk process is defined by
\BQNY
R(t)=x+pt -\sigma B(t),\ \ t\ge 0,
\EQNY
where $x>0$ is the {\it initial capital}, $p>0$ is the {\it net profit rate} and $\sigma B(t)$ models the net loss process with    $\sigma>0$ the volatility coefficient. Roughly speaking,  $\sigma B(t)$ is an approximation of the total claim amount process by time $t$ minus its expectation, the latter is usually called the {\it pure premium} amount and calculated  to cover the   average payments of claims. The net profit, also called {\it safety loading}, is the component which protects the company from large deviations of claims from the average and also allows an accumulation of capital.
Ruin related problems for  Brownian models  have been well studied; see, e.g., \cite{GS04, AsmAlb10}.


In recent years, multi-dimensional risk models have been introduced to model the surplus of multiple business lines  of an insurance company or the suplus of collaborating companies (e.g., insurance and reinsurance). We refer to \cite{AsmAlb10} [Chapter XIII 9] and
\cite{AL18, AM17, APP08a, APP08b, FKPR17, JR18, AMP16, AAM17,  AM18} for relevant  recent discussions. It is concluded in the literature that in comparison with the well-understood 1-dimensional risk models, study of multi-dimensional risk models is much more challenging.  It is shown recently in \cite{DMSW18} that multi-dimensional Brownian model can serve as approximation of a multi-dimensional classical risk model in a Markovian environment. Therefore, obtained results for multi-dimensional Brownian model can serve as approximations of  the multi-dimensional  classical risk models in a Markovian environment; ruin probability approximation has been used in the aforementioned paper.
Actually, multi-dimensional Brownian models have drawn a lot of attention due to its tractability and practical relevancy.

A   $d$-dimensional  Brownian  model can be defined in a matrix form as
\BQNY
\vk R(t) = \vk x +\vk p t - \vk X(t),  \ \ t\ge 0, \ \ \text{with}\ \ \vk X(t)=A\vk{B}(t),
\EQNY
where $\vk x=(x_1,\cdots,x_d)^\top, \vk p=(p_1,\cdots,p_d)^\top\in (0,\IF)^d$ are, respectively, (column) vectors representing the initial capital and net profit rate,
$A \in \R^{d\times d}$ is a non-singular matrix modelling dependence between different business lines, and
$\vk{B}(t)=(B_1(t),\ldots,B_d(t))^\top,t\ge0$ is a standard $d$-dimensional Brownian motion (BM)
with independent coordinates. Here $\top$ is  the transpose sign. In what follows, vectors are understood as column vectors written in bold letters. 

Different types of ruin can be considered in multi-dimensional models, which are relevant to the probability that the
surplus of one or more of the business lines drops below zero in a certain time interval $[0,T]$ with $T$ {either}
a finite constant or infinity. \kk{One of the}  commonly 
studied is the so-called {\it simultaneous ruin probability} defined as
\[
Q_T(\vk x):=
\pk { \exists_{t\in[0,T]}\bigcap_{i=1}^d \Bigl\{ R_i(t)<0\Bigl\}},
\]
which is the probability that at a certain time $t\in[0,T]$ all the surpluses become negative.
Here for $T<\IF$, $Q_T(\vk x)$ is called finite-time \kk{simultaneous} ruin probability,
and $Q_\IF(\vk x)$ is called infinite-time \kk{simultaneous} ruin probability. Simultaneous ruin probability, which is essentially the hitting probability of {$\vk R (t)$ to the orthant $\{\vk{y}\in\R^d: y_i<0, i=1,\ldots,d\}$}, has been discussed for multi-dimensional Brownian models in
different contexts; see \cite{Garbit14, DHJT18}.  In \cite{Garbit14}, for  fixed $\vk x$ the asymptotic behaviour of $Q_T(\vk x)$ as $T\to\IF$ has been discussed. Whereas, in \cite{DHJT18}, the asymptotic behaviour, as $u\to\IF$, of the infinite-time ruin probability $Q_\IF(\vk x)$, with $\vk{x}=\vk \alpha u=(\alpha_1 u, \alpha_2 u, \ldots, \alpha_d u)^\top, \alpha_i>0, 1\le i\le d$
has been obtained. Note that it is common  in risk theory to derive the later  type of asymptotic results for  ruin probabilities; see, e.g., \cite{EKM97, Mik08, APP08a}.

Another type of ruin probability is the {\it component-wise (or joint) ruin probability} defined as
\BQN\label{eq:PPu0}
P_T(\vk x):=\pk {\bigcap_{i=1}^d \Bigl\{ \exists_{t\in[0,T]} R_i(t)<0 \Bigl\}}
=\pk {\bigcap_{i=1}^d \Bigl\{ \sup_{t_i\in[0,T]} (X_i(t_i)-p_i t_i)>x_i \Bigl\}},
\EQN
which is the probability that all surpluses get below zero but possibly at different times.  
It is  this possibility that makes the study of $P_T(\vk x)$ more difficult.

\kk{The study} of  joint distribution of the extrema of multi-dimensional BM over finite-time interval
has been proved to be important in quantitative finance; see, e.g.,  \cite{HKR98, KZ16}.
We refer to  \cite{DMSW18} for a comprehensive summary of related results. Due to the complexity of the problem,
two-dimensional case has been the focus in the literature, and for this case some explicit \kk{formulas} can be obtained
by using a PDE approach. Of
particular relevance to 
\kk{the} ruin probability $P_T(\vk x)$ is a result derived in
\cite{HKR98} which shows that
\BQNY
&&\pk{\sup_{t\in[0,T]} ( X_1(t) - p_1 t ) \le  x_1, \ \sup_{s\in[0,T]}  (X_2(s) - p_2 s) \le  x_2 }\\
&&=e^{a_1 x_1 +a_2 x_2+b  T} f(x_1,x_2, T),
\EQNY
where $a_1,a_2,b$  are known constants, and $f$ is a function defined in terms of infinite-series, double-integral and Bessel function.  Using the above formula one can derive an expression for $P_T(\vk x)$  in two-dimensional case as follows
\BQN \label{eq:QT}
P_T(\vk x)&=&1-\pk{\sup_{t\in[0,T]} ( X_1(t) - p_1 t ) \le  x_1}-\pk{\sup_{s\in[0,T]}  (X_2(s) - p_2 s) \le  x_2 } \\
&& +\pk{\sup_{t\in[0,T]} ( X_1(t) - p_1 t ) \le  x_1, \ \sup_{s\in[0,T]}  (X_2(s) - p_2 s) \le  x_2 },\nonumber
\EQN
where the expression for the distribution of single supremum is also known; see \cite{HKR98}.
 Note that even though  we have obtained explicit expression of $P_T(\vk x)$ in \eqref{eq:QT} for the two-dimensional case,  it seems \kk{difficult to derive the explicit form of} the corresponding  infinite-time ruin probability $P_\IF(\vk x)$ by simply putting $T\to\IF$ in \eqref{eq:QT}.

By assuming $\vk{x}=\vk \alpha u=(\alpha_1 u, \alpha_2 u, \ldots, \alpha_d u)^\top, \alpha_i>0, 1\le i\le d$, we aim to analyse the asymptotic behaviour of the infinite-time ruin probability $P_\IF(\vk x)$ as $u\to\IF$.  Applying
Theorem 1 in \cite{RolskiSPA} we \kk{arrive at} the following logarithmic asymptotics 
\BQN \label{eq1D}
-\frac{1}{u}\ln P_\IF(\vk x) &\sim &  \frac{1}{2}\inf_{\vk t > \vk 0}
  \inf_{\vk{v} \ge \vk{\alpha}+\vk{p} \vk t}  \vk{v}^\top \Sigma^{-1}_{\vk t}  \vk{v}, \ \ \ \ \text{as} \ u\to\IF
\EQN
 provided $\Sigma_{\vk t}$ is non-singular,  where $\vk{pt}:=(p_1 t_1,\cdots, p_d t_d)^\top$, inequality of vectors are meant component-wise,
 and $\Sigma^{-1}_{\vk t}$ is the inverse matrix of  the covariance function $ \Sigma_{\vk t}$ of $(X_1(t_1), \cdots, X_d(t_d))$, with $\vk t =(t_1,\cdots, t_d)^\top$ and $\vk 0 =(0,\cdots, 0)^\top\in \R^d$. Let us recall that conventionally
for two given positive functions $f(\cdot)$ and $h(\cdot)$, we write
 $f(x)\sim h(x)$ if  $ \lim_{x \to \IF}  {f(x)}/{h(x)} = 1$. 

 For {more precise analysis on $P_\IF(\vk x)$}, it seems crucial to first solve the  two-layer optimization problem in \eqref{eq1D} and
 find the optimization points $\vk t_0$.  As it can be {recognized} in the following, when dealing with $d$-dimensional
 case with $d>2$ \kk{the calculations become highly nontrivial and} complicated.
 Therefore, in this contribution we only discuss a tractable two-dimensional model and aim for an explicit logarithmic
 asymptotics by solving the minimization problem in \eqref{eq1D}.

 In the classical ruin theory when analysing the compound Poisson model or Sparre Andersen model,
the so-called {\it adjustment coefficient}  is used as a measure of goodness; see, e.g., \cite{AsmAlb10} or \cite{Rolskibook}.   It is of interest to obtain the solution of the minimization problem in \eqref{eq1D} from a practical point of view, as {it can be seen as an analogue of the adjustment coefficient and thus we could get some insights about the risk that the company is facing.}
As discussed in \cite{AsmAlb10} and \cite{LLT07} it is also of interest to know how the dependence between different risks influences the joint ruin probability,  which can be easily analysed through the obtained logarithmic asymptotics; see Remark \ref{Rem1}.


  The rest of this paper is organised as follows. In
 Section \ref{MM},
we formulate the two-dimensional  Brownian model and give the main results of this paper.
The main lines of proof with  auxiliary lemmas are displayed in  Section \ref{Sec:FRP}.  In Section \ref{Sec:CD} we conclude the paper.
All technical proofs of the lemmas in Section \ref{Sec:FRP} are presented in  Appendix.


 \section{Model formulation and main results} \label{MM}
\Lcc{Due to the fact that component-wise ruin probability $P_\IF(\vk x)$ does not change under scaling, we can simply assume that the volatility coefficient for all business lines is equal to 1. Furthermore, noting that the timelines for different business lines should be distinguished as shown in \eqref{eq:PPu0} and \eqref{eq1D}, we  introduce} a two-parameter extension of correlated two-dimensional BM defined as
\[\left(X_1(t),X_2(s)\right)=\left(B_1(t),\ \rho B_1(s)+\sqrt{1-\rho^2}B_2(s)\right), \ \ t,s\ge0,
\]
with $\rho\in
(-1,1)$ and mutually independent Brownian motions $B_1,B_2$.  We shall consider the following two dependent insurance risk processes
$$
R_i(t)=u+\mu_i t-X_i(t), \ \ t\ge0,\ \ \ \ i=1,2,
$$
where $\mu_1,\mu_2>0$ \Lcc{are net profit rates, $u$ is the  initial capital (which is assumed to be the same for both business  lines, as otherwise, the calculations become rather complicated).}
We shall assume without loss of generality that $\mu_1\le \mu_2$. \Lcc{Here, $\mu_i$ is different from $p_i$
(see (\ref{eq:PPu0})) in  the sense that it corresponds to the (scaled) model with volatility coefficient
 standardized to be 1.}


In this contribution,  we shall focus  on the logarithmic
asymptotics of
\Lcc{\BQN\label{eq:PPu1}
P(u):=P_\IF(u(1,1)^\top)
&=&
\pk{\{\exists_{t\ge 0} R_1(t)<0\}\ \cap \ \{\exists_{s\ge 0} R_2(s)<0\}}\\
&=&
\pk{\sup_{t\ge 0} ( X_1(t) - \mu_1 t ) >  u, \ \sup_{s\ge 0}  (X_2(s) - \mu_2 s) > u }, \  \ {\rm as} \ u\to\IF. \nonumber
\EQN
}




Define
\BQN\label{eq:rho12}
\hat\rho_1=\frac{\mu_1+\mu_2-\sqrt{(\mu_1+\mu_2)^2-4\mu_1(\mu_2-\mu_1)}}{4\mu_1}\in [0,\frac{1}{2}), \ \ \ \
\hat\rho_2=\frac{\mu_1+\mu_2}{2\mu_2}
\EQN
and let
\BQN\label{eq:tsstar}
 t^*=t^*(\rho)= s^*= s^*(\rho):=\sqrt {\frac{2(1-\rho)}{\mu_1^2+\mu_2^2-2\rho\mu_1\mu_2}}.
\EQN

The following theorem constitutes the main result of this contribution.

\begin{theo}\label{p.log}  \jc{For the joint infinite-time ruin probability \eqref{eq:PPu1} we have, as $u\to\IF$},
\begin{eqnarray*}
-\frac{\log(P(u))}{u}\sim
\left\{
  \begin{array}{ll}
    2(\mu_2+(1-2\rho)\mu_1), & \hbox{if \ $-1< \rho  \le\hat \rho_1$;} \\
    \frac{\mu_1+\mu_2+2/t^*}{1+\rho}, & \hbox{if \ $\hat \rho_1  < \rho<\hat \rho_2$ ;} \\
    2\mu_2, & \hbox{if \ $\hat \rho_2\le\rho<1$.}
  \end{array}
\right.
\end{eqnarray*}
\end{theo}

\begin{remark} \label{Rem1}
\Lcc{ a).  Following the classical one-dimensional risk theory we can call
   quantities on the right hand side in Theorem \ref{p.log} as {\it adjustment coefficients}.
   They serve sometimes as a measure of goodness for a risk business.}
   
 \Lcc{b). One can easily check that    adjustment coefficient as a function of $\rho$
 is continuous,  strictly decreasing  on $(-1,\hat \rho_2]$, and it is constant,  {equal to} $2\mu_2$ on $[\hat \rho_2, 1)$. This means that as the two lines of business becomes more positively correlated the risk of ruin becomes larger, which is  consistent with the intuition.
}\end{remark}

Define
\BQN \label{eq:gqpp}
g(t,s):= \underset{y\ge 1+\mu_2 s}{\inf_{x\ge 1+\mu_1 t}}\  (x,y)\ \Sigma_{ts}^{-1}\  (x,y)^\top,\ \ \ \ \ t,s>0,
\EQN
where $\Sigma_{ ts}^{-1}$ is the inverse matrix of
$
\  \Sigma_{ts}=
\left(
   \begin{array}{cc}
     t & \rho\  t\wedge s \\
     \rho\  t\wedge s & s \\
   \end{array}
 \right),
$
with $ t\wedge s=\min(t,s)$ and { $\rho\in(-1,1)$}.

The proof of Theorem \ref{p.log} \Lcc{follows from \eqref{eq1D}} 
which implies that
the logarithmic asymptotics for $P(u)$ is of the form
\BQN\label{eq:logri1}
-\frac{1}{u}\ln P(u) &\sim& \frac{g(\vk{t}_0)}{2}, \ \ \ u\to\IF,
\EQN
where
\BQN\label{eq:ght}
g(\vk{t}_0)=\inf_{(t,s)\in(0,\IF)^2} g(t,s),
\EQN
and Proposition \ref{Lem:Optg} below, wherein we list dominating points $\vk t_0$
that optimize  \Lcc{the function $g$ over $(0,\IF)^2$ and the corresponding optimal values $g(\vk t_0)$.}

In order to solve the \jc{two-layer minimization problem} in \eqref{eq:ght} (see also \eqref{eq:gqpp})
 we define for $t,s>0$ the following functions:
\BQNY
&&g_1(t)=\frac{(1+\mu_1 t)^2}{t}, \ \ \ \ g_2(s)=\frac{(1+\mu_2 s)^2}{s},\\
 && g_3(t,s)=(1+\mu_1 t, 1+\mu_2 s) \ \Sigma_{ts}^{-1} \  (1+\mu_1 t, 1+\mu_2 s)^\top.
\EQNY
Since $t \wedge s$ appears in the above formula, we shall consider a partition of the quadrant $(0,\IF)^2$, namely
\BQN\label{eq:AB}
(0,\IF)^2=A\cup L\cup B,\ \ \ A=\{s< t\},\ L=\{s=t\},\  B=\{s>t\}.
\EQN
For convenience we denote $ \oA=\{s\le t\}=A\cup L$ and  $ \oB=\{s\ge t\}=B\cup L$. Hereafter, all  sets are defined on $(0,\IF)^2$, so $(t,s)\in (0,\IF)^2$ will be omitted.


Note that $g_3(t,s)$ can be represented in the following  form:
\BQN
g_3(t,s)
\COM{&=&\left\{\begin{array}{cc}
 \frac{ (1+\mu_1 t)^2 s -2\rho s  (1+\mu_1 t)(1+\mu_2 s) + (1+\mu_2s)^2 t }{ts-\rho^2 s^2}, & \hbox{if }  (t,s)\in \oA \label{eq:g31}\\
\frac{ (1+\mu_1 t)^2 s -2\rho t  (1+\mu_1 t)(1+\mu_2 s) + (1+\mu_2s)^2 t }{ts-\rho^2 t^2}, , & \hbox{if }   (t,s)\in \oB
\end{array}\right.\\}
&=&\left\{\begin{array}{cc}
g_A(t,s):=\frac{(1+\mu_2s)^2}{s} +\frac{((1+\mu_1 t)-\rho(1+\mu_2 s))^2}{t-\rho^2 s}, & \hbox{if }  (t,s)\in \oA \\
g_B(t,s):=\frac{(1+\mu_1 t)^2}{t} +\frac{((1+\mu_2 s)-\rho(1+\mu_1 t))^2}{s-\rho^2 t}, & \hbox{if }   (t,s)\in \oB.
\end{array}\right.\label{eq:g32}
\EQN
Denote further
\BQN\label{eq:fL}
g_L(s):=g_A(s,s)=g_B(s,s)=\frac{(1+\mu_1s)^2 +(1+\mu_2s)^2-2\rho(1+\mu_1s)(1+\mu_2s) }{(1-\rho^2) s},\ \ \ s>0.
\EQN

In the next proposition we identify the so-called dominating points,
that is, points \jc{$\vk t_0$} for which function defined in \eqref{eq:gqpp} achieves its minimum. This identification
 might also be  useful for deriving a more subtle asymptotics for $P(u)$.

{\bf Notation:}
{\it In the following, in order to keep the notation consistent, $\rho\le \mu_1/\mu_2$ is understood as $\rho<1$ if $\mu_1=\mu_2.$}


\begin{prop}
 \label{Lem:Optg} 
\begin{itemize}
\item[(i).] Suppose that  $-1< \rho < 0$. \\
For    $\mu_1<\mu_2$ we have
\BQNY
 g(\vk{t}_0)  = g_A(t_A,s_A) =   4(\mu_2+(1-2\rho)\mu_1),
\EQNY
where,
$
(t_A,s_A)=(t_A(\rho),s_A(\rho)):=\LT(\frac{1-2\rho}{\mu_1}, \frac{1}{\mu_2-2\mu_1\rho} \RT)
$
is the unique minimizer of $g(t,s), (t,s)\in(0,\IF)^2$.
\\For $\mu_1=\mu_2=:\mu$ we have
\BQNY
 g(\vk{t}_0) = g_A(t_A,s_A) = g_B(t_B,s_B)=  \jc{8(1- \rho)\mu,}
\EQNY
where
$
(t_A,s_A)=\LT(\frac{1-2\rho}{\mu}, \frac{1}{ (1-2 \rho)\mu} \RT)\in A,
 (t_B,s_B):=\LT(\frac{1}{ (1-2 \rho)\mu},\frac{1-2\rho}{\mu} \RT)\in B$
are the only two  minimizers of $g(t,s), (t,s)\in(0,\IF)^2$.

\item[(ii).] Suppose that $0\le \rho<  \hat \rho_1$. We have
\BQNY
g(\vk{t}_0)  = g_A(t_A,s_A) =   4(\mu_2+(1-2\rho)\mu_1),
\EQNY
where $(t_A,s_A)\in A$ is the unique minimizer  of $g(t,s), (t,s)\in(0,\IF)^2$.  

\item[(iii). ] Suppose that $\rho = \hat \rho_1$. We have
\BQNY
g(\vk{t}_0) = g_A(t_A,s_A) =   4(\mu_2+(1-2\rho)\mu_1),
\EQNY
where $(t_A,s_A) =(t_A(\hat\rho_1), s_A(\hat\rho_1))=(t^*(\hat\rho_1), s^*(\hat\rho_1))\in L$,
is the unique minimizer  of $g(t,s), (t,s)\in(0,\IF)^2$, with $(t^*,s^*)$ defined in \eqref{eq:tsstar}.

\item[(iv).] Suppose that $ \hat \rho_1  < \rho< \hat \rho_2 $. We have
\BQNY
g(\vk{t}_0)   = g_A(t^*, s^{*}) = g_L(t^*)= \frac{2}{1+\rho}(\mu_1+\mu_2+2/t^*),
\EQNY
where $(t^*,s^*)\in L$ is  the unique minimizer  of $g(t,s), (t,s)\in(0,\IF)^2$.

\item[(v).] Suppose that $ \rho  =  \hat \rho_2 $. We have    $t^*(\hat\rho_2)=s^*(\hat\rho_2)=1/\mu_2$, and
\BQNY
g(\vk{t}_0) = g_A(1/\mu_2 ,1/\mu_2 )=g_L( 1/\mu_2 ) =g_2( 1/\mu_2 )=4\mu_2,
\EQNY
where
the  minimum  of $g(t,s), (t,s)\in(0,\IF)^2$ is attained at $(1/\mu_2 ,1/\mu_2 )$,  with
$g_3(1/\mu_2 ,1/\mu_2 )=g_2(1/\mu_2)$, and $1/\mu_2 $ is the unique minimizer of $g_2(s), s\in(0,\IF)$.

\item[(vi).] Suppose that  $ \hat \rho_2<\rho <1$. We have
\BQNY
g(\vk{t}_0)  = 
g_2( {1}/{\mu_2})= 4\mu_2,
\EQNY
\jc{where
the  minimum  of $g(t,s), (t,s)\in(0,\IF)^2$ is attained when $g(t,s)=g_2(s)$.}
\end{itemize}
\end{prop}





 \begin{remark}\label{Rem:gg} 
 In case that $\mu_1=\mu_2$, we have $\hat\rho_1=0, \hat\rho_2=1$, and thus
 scenarios (ii)  and  (vi) do not apply.
%
%
 \end{remark}


\section{Proofs of main results} \label{Sec:FRP}

\jc{As discussed in the previous section, Proposition \ref{Lem:Optg} combined with (\ref{eq:logri1}), straightforwardly implies the thesis of Theorem \ref{p.log}. In what follows, we shall focus on the proof of Proposition \ref{Lem:Optg}, for which we need to find the dominating points $\vk t_0$ by solving the two-layer minimization problem (\ref{eq:ght}).}

The solution of quadratic programming problem of the form \eqref{eq:gqpp} \jc{(inner minimization problem of \eqref{eq:ght})}  has been well understood; e.g., \cite{HA2005, ENJH02} (see also Lemma 2.1 of \cite{DHJT18}). For completeness and for reference, we present below Lemma 2.1 of \cite{DHJT18} for the case where $d=2$.

We introduce some more notation. If $I\subset \{1,2\}$, then for a vector $\vk a\in\R^2$ we denote by  $\vk{a}_I=(a_i, i\in I)$ a sub-block vector of $\vk a$. Similarly, if further $J \subset \{1,2\}$, for  a matrix $M=(m_{ij})_{i,j\in \{1,2\}}\in \R^{2\times 2}$ we denote  by  $M_{IJ}\ccj{=M_{I,J}}=(m_{ij})_{i\in I, j\in J}$  the  sub-block matrix of $M$
 determined by $I$ and $J$.  Further, write $M_{II}^{-1}=(M_{II})^{-1}$ for the inverse matrix of $M_{II}$ whenever it exists.

\BEL \label{AL}
Let
$M \in \R^{2 \times 2}$ be a positive definite  matrix. 
If  $\vk{b}\in  \R ^2 \setminus (-\infty, 0]^2 $, then the quadratic programming problem
$$ P_M(\vk{b}): \text{Minimise $ \vk{x}^\top M^{-1} \vk{x} $ under the linear constraint } \vk{x} \ge \vk{b} $$
has a unique solution $\widetilde{\vk{b}}$ and there exists a unique non-empty
index set $I\subseteq \{1, 2\}$ such that
\begin{eqnarray*}  \label{eq:IJi}
&&\widetilde{\vk{b}}_{I} =
\vk{b}_{I}\not= \vk 0_I, \quad M_{II}^{-1} \vk{b}_{I}>\vk{0}_I,\\
\text{and}&& \text{if}\  I^c =\{ 1, 2\} \setminus I \not=
\emptyset, \text{ then }
\widetilde {\vk{b}}_{I^c}
=  M_{I^cI}M_{II}^{-1} \vk{b}_{I}\ge \vk{b}_{I^c}.
\label{eq:hii}
\end{eqnarray*}
Furthermore,
\BQNY
\min_{\vk{x} \ge  \vk{b}}
\vk{x}^\top M^{-1}\vk{x}&=& \widetilde{\vk{b}}^\top M^{-1} \widetilde {\vk{b}}   =
\vk{b}_{I}^\top M_{II}^{-1}\vk{b}_{I}>0, \\
 \label{eq:new}
\vk{x}^\top M^{-1} \widetilde{ \vk{b}}&=& \vk{x}_I^\top M_{II}^{-1} \widetilde {\vk{b}}_I=
\vk{x}_I^\top M_{II}^{-1}\vk{b}_I,  \quad
\forall \vk{x}\in \R^2.
\EQNY
\EEL

{For the solution of the quadratic programming problem \eqref{eq:gqpp} a} suitable representation for $g(t,s)$ is worked out in the following lemma.

\Lcc{For $1>\rho> \mu_1/\mu_2,$} let $D_2=\{(t,s): w_1(s)\le t\le f_1(s)\}$ and $D_1=(0,\IF)^2\setminus D_2$ ,
with boundary functions given by
\BQN\label{eq:fw}
f_1(s)=\frac{\rho-1}{\mu_1}+\frac{\rho \mu_2}{\mu_1}s, \ \  w_1(s)= \frac{s}{\rho+( \rho \mu_2-\mu_1)s}, \ \ s\ge0,
\EQN
and the unique intersection point of $f_1(s), w_1(s), s\ge 0,$   given by
\BQN\label{eq:s1}
s_1^* =s_1^*(\rho) : =\frac{1-\rho}{\rho \mu_2 -\mu_1},
\EQN
as depicted  in Figure 1.

 \begin{figure}
 \vspace{0mm}
  \includegraphics[width=100mm, height=60mm]{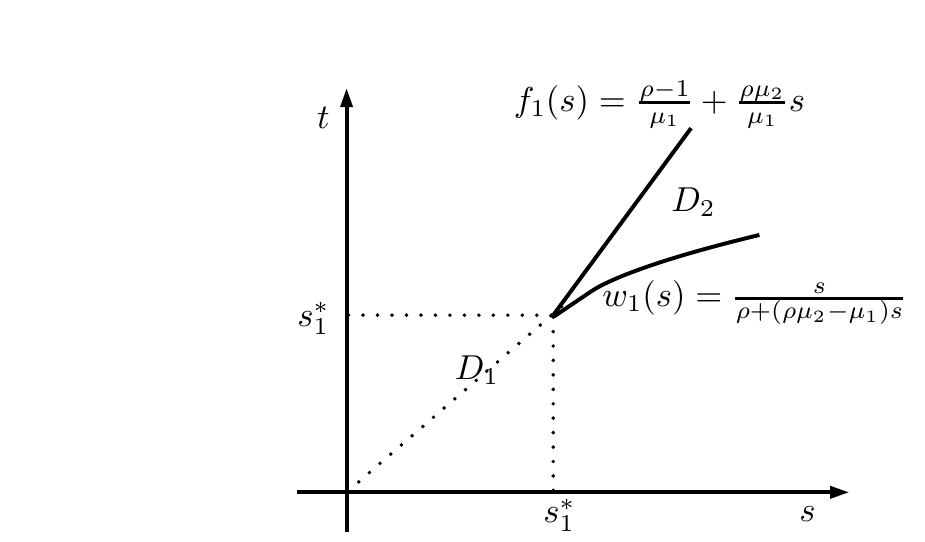}
  \caption{Partition of $(0,\IF)^2$ into $D_1, D_2$}
\end{figure}

\BEL \label{Lem:qpp}
Let $g(t,s), t,s>0$ be given as in \eqref{eq:gqpp}. We have:
\begin{itemize}
\item[(i).] If $-1< \rho\le \mu_1/\mu_2,$ then
$$
g(t,s)=g_3(t,s),\ \ \ (t,s)\in (0,\IF)^2.
$$
\item[(ii).] If  $1>\rho> \mu_1/\mu_2,$ then
$$
g(t,s)=\left\{\begin{array}{cc}
g_3(t,s), & \hbox{if }  (t,s)\in D_1\\
 g_2(s), & \hbox{if }   (t,s)\in D_2.
\end{array}\right.
$$
\end{itemize}

Moreover, we have $g_3(f_1(s), s)=g_3(w_1(s), s)= g_2(s)$ for all $s\ge s_1^*$.
\EEL
\subsection{Proof of Proposition \ref{Lem:Optg}}
We shall discuss in order the case when \jc{$-1<  \rho<0$} 
and the case when $0\le  \rho<1$ in the following two subsections.
In both scenarios we shall  first derive the minimizers of the function $g(t,s)$ on regions $\oA$ and $\oB$ (see \eqref{eq:AB}) separately, and then look for a global minimizer by comparing the two minimum values.
\jc{For clarity some scenarios}  are analysed in forms of lemmas.
\subsubsection{Case $-1<  \rho<0$} By Lemma \ref{Lem:qpp}, we have that
\BQNY
g(t,s)=g_3(t,s), \ \ \  (t,s)\in (0,\IF)^2.
\EQNY
We shall derive the minimizers of $g_3(t,s)$ on $\oA, \oB$ separately.

\underline{Minimizers of $g_3(t,s)$ on  $\oA$}. We have,  for any fixed $s$,
\BQNY
\frac{\partial g_3(t,s)}{\partial t}=\frac{\partial g_A(t,s)}{\partial t}=0 \ \ \Leftrightarrow \ \ (\mu_1t+1-\rho-\rho\mu_2 s)(\mu_1 t -(2\mu_1\rho^2-\rho \mu_2)s+\rho-1)=0,
\EQNY
where the representation \eqref{eq:g32} is used.
Two roots of the above equation are:
\BQN\label{eq:tt}
t_1=t_1(s):=\frac{\rho-1+\rho\mu_2s}{\mu_1},\ \ \ \ \ t_2=t_2(s):=\frac{1-\rho+(2\mu_1\rho^2-\rho \mu_2)s}{\mu_1}.
\EQN
Note that, due to the form of the function $g_A(t,s)$   given in \eqref{eq:g32}, for any fixed $s$, there exists a unique minimizer of $g_A(t,s)$ on $\oA$ which is either an inner point $t_1$ or $t_2$ (the one that is larger than $s$), or a boundary point $s$.
Next, we check if any of $t_i, i=1,2,$ is larger than $s$. 
Since $\rho<0$, $t_1<0<t_2$. So we check if $t_2>s$. It can be shown that
\BQN\label{eq:t2s}
t_2>s \ \ \Leftrightarrow \ \  (\mu_1+\rho \mu_2-2\mu_1\rho^2)s<1-\rho.
\EQN
Two scenarios $\mu_1+\rho \mu_2-2\mu_1\rho^2\le 0$ and $\mu_1+\rho \mu_2-2\mu_1\rho^2> 0$ will be distinguished.

\underline{Scenario $\mu_1+\rho \mu_2-2\mu_1\rho^2\le 0$}.   We have from \eqref{eq:t2s} that
$$
t_1<0<s<t_2,
$$
and thus
$$
\inf_{(t,s)\in \oA} g_3(t,s) =\inf_{s>0} f_A(s),
$$
where
\BQNY
f_A(s):=g_A(t_2(s),s)=\frac{(1+\mu_2 s)^2}{s}+4\mu_1((1-\rho)+(\rho^2\mu_1-\rho\mu_2)s).
\EQNY
Next, since
\BQN\label{eq:sA}
f_A'(s)=0 \ \  \Leftrightarrow \ \ s_A=s_A(\rho):=\frac{1}{\abs{\mu_2-2\rho\mu_1}}=\frac{1}{\mu_2-2\rho\mu_1}>0,
\EQN
the unique minimizer of $g_3(t,s)$ on $\oA$ is given by $(t_A, s_A)$ with
$$
 t_A:=t_2(s_A)=\frac{1-2\rho}{\mu_1}.
$$

\underline{Scenario $\mu_1+\rho \mu_2-2\mu_1\rho^2> 0$}.  We have from \eqref{eq:t2s} that
\BQN\label{eq:SDoubs}
t_1<0<s<t_2\ \ \Leftrightarrow \ \ s<\frac{1-\rho}{\mu_1+\rho \mu_2-2\mu_1\rho^2}=\frac{1-\rho}{\rho(\mu_2-\mu_1\rho)+\mu_1(1-\rho^2) }=:s^{**}(\rho)=s^{**}, %
\EQN
and in this case,
\BQN\label{eq:g3A}
\inf_{(t,s)\in \oA} g_3(t,s) =\min\LT(\inf_{0<s<s^{**}} f_A(s),  \inf_{s\ge s^{**}} g_L(s)\RT),
\EQN
where $g_L(s)$ is given in \eqref{eq:fL}.
Note that
\BQN\label{eq:stStar}
g_L'(s)=0 \ \  \Leftrightarrow \ \ s^*=s^*(\rho)= \sqrt {\frac{2(1-\rho)}{\mu_1^2+\mu_2^2-2\rho\mu_1\mu_2}}.
\EQN
Next, 
\Lcc{for $-1< \rho<0$ we have} that (recall $s^{**}$ given in \eqref{eq:SDoubs}) 
\BQNY
s^{**}\Lcc{\ge \frac{1-\rho}{\mu_1(1-\rho^2)}}>\frac{1}{\mu_1}\ge \frac{1}{\mu_2}>s_A, \ \ \
s^{**}>\frac{1-\rho}{\mu_1} >s^*.
\EQNY
Therefore, by  \eqref{eq:g3A} we conclude that  the unique minimizer of $g_3(t,s)$ on $\oA$ is again given by $(t_A, s_A)$.
Consequently, for all $-1< \rho<0$,  we have that  the unique minimizer of $g_3(t,s)$ on $\oA$ is given by $(t_A, s_A)$, and
\BQN\label{eq:Ag}
\inf_{(t,s)\in \oA} g_3(t,s) =  g_A(t_A,s_A)=4(\mu_2+(1-2\rho)\mu_1).
\EQN

{
\underline{Minimizers of $g_3(t,s)$ on $\oB$}.
Similarly, we have, for any fixed $t$,
\BQNY
\frac{\partial g_3(t,s)}{\partial s}=\frac{\partial g_B(t,s)}{\partial s}=0 \ \ \Leftrightarrow \ \ (\mu_2 s+1-\rho-\rho\mu_1 t)(\mu_2s -(2\mu_2\rho^2-\rho \mu_1)t+\rho-1)=0.
\EQNY
Two roots of the above equation are:
\BQN \label{eq:ss}
s_1=s_1(t):=\frac{\rho-1+\rho\mu_1 t}{\mu_2},\ \ \ \ \ s_2=s_2(t):=\frac{1-\rho+(2\mu_2\rho^2-\rho \mu_1)t}{\mu_2}.
\EQN
Next, we check if any of $s_i, i=1,2,$ is greater than $t$. 
Again $s_1<0<s_2$ as $\rho<0$. So we check if $s_2>t$. It can be shown that
\BQN\label{eq:t2t}
s_2>t \ \ \Leftrightarrow \ \  (\mu_2+\rho \mu_1-2\mu_2\rho^2) t<1-\rho.
\EQN
Thus, for \underline{Scenario $\mu_2+\rho \mu_1-2\mu_2\rho^2\le 0$} we have that
$$
s_1<0<t<s_2
$$
and in this case
$$
\inf_{(t,s)\in \oB} g_3(t,s) =\inf_{t>0} f_B(t),
$$
with
\BQNY 
f_B(t):=g_B(t,s_2(t))=\frac{(1+\mu_1 t)^2}{t}+4\mu_2((1-\rho)+(\rho^2\mu_2-\rho\mu_1)t).
\EQNY
Next, note that
\BQN\label{eq:tB}
f_B'(t)=0 \ \  \Leftrightarrow \ \ t_B=t_B(\rho):=\frac{1}{\abs{\mu_1-2\rho\mu_2}}=\frac{1}{\mu_1-2\rho\mu_2}>0.
\EQN
Therefore,   the unique minimizer of $g_3(t,s)$ on $\oB$ is given by $(t_B, s_B)$ with
\BQNY
 s_B:=s_2(t_B)=\frac{1-2\rho}{\mu_2}, \ \ \ \ \inf_{(t,s)\in \oB} g_3(t,s) =  g_B(t_B,s_B)=4(\mu_1+(1-2\rho)\mu_2).
\EQNY

For \underline{ Scenario $\mu_2+\rho \mu_1-2\mu_2\rho^2> 0$} we have from \eqref{eq:t2t} that
\BQN\label{eq:TDoubs}
s_1<0<t<s_2\ \ \Leftrightarrow \ \ t<\frac{1-\rho}{\mu_2+\rho \mu_1-2\mu_2\rho^2}=\frac{1-\rho}{\rho(\mu_1-\rho \mu_2)+\mu_2(1-\rho^2)}=:t^{**}(\rho)=t^{**}. 
\EQN
In this case,
\BQNY
\inf_{(t,s)\in \oB} g_3(t,s) =\min\LT(\inf_{0<t<t^{**}} f_B(t),  \inf_{t\ge t^{**}} g_L(t)\RT).
\EQNY
Though it is not easy to determine explicitly the optimizer, we can conclude that
the minimizer should be taken at  $(t_B,s_B)$, $(t^*,t^*)$ or $(t^{**},t^{**})$, where $t^*=t^*(\rho)=s^*(\rho)$. Further, we have
 from the discussion in \eqref{eq:g3A} that
$$
g_A(t_A,s_A)<g_L(s^*)=g_L(t^*)=\min(g_L(t^*),g_L(t^{**})),
$$
and
$$
g_B(t_B,s_B)=4(\mu_1+(1-2\rho)\mu_2)\ge 4(\mu_2+(1-2\rho)\mu_1) = g_A(t_A,s_A).
$$
Combining the above discussions on $\oA, \oB$, we conclude that Proposition
\ref{Lem:Optg} holds for $-1< \rho<0$.
}\QED

{
\subsubsection{Case $0\le \rho<1$}

We shall derive the minimizers of $g(t,s)$ on $\oA,\oB$ separately. We start with discussions on  $\oB$, for which we give the following lemma. Recall $t^*(\rho)=s^*(\rho)$ defined in \eqref{eq:stStar} (see also \eqref{eq:tsstar}), $t_B(\rho)$ defined in \eqref{eq:tB}, $t^{**}(\rho)$ defined in \eqref{eq:TDoubs}, and $s_1^*(\rho)$ defined in \eqref{eq:s1} for $\mu_1/\mu_2<\rho<1$.
Note that where it applies,  $1/0$ is understood as $+\IF$ and $1/\IF$ is understood as $0$.



\BEL\label{Lem:tt}
We have:

(a).    The function $t^*(\rho)$ is a decreasing function  on $[0, 1]$, and both $t_B(\rho)$ and $s_1^*(\rho)$ are  decreasing functions on $(\mu_1/\mu_2,1)$.

(b).  The function  $t^{**}(\rho)$ decreases from $1/\mu_2$ at $\rho=0$ to some positive value and then increases to  $1/\mu_2$ at $\hat\rho_2$ (defined in \eqref{eq:rho12}), and then increases to $+\IF$ at the root $\Lcc{\hat\rho}\in(0,1]$ of  the equation $\mu_2+\rho\mu_1-2\mu_2\rho^2=0.$

(c). For $0\le \rho\le \mu_1/\mu_2$, we have
\BQNY 
t_B(\rho)\ge t^{**}(\rho), \ \ \ t^*(\rho)\ge t^{**}(\rho),
\EQNY
where both equalities hold only when $\rho=0$ and $\mu_1=\mu_2.$

(d). It holds that
\BQN\label{eq:tttt}
t^*(\hat\rho_2)=t_B(\hat\rho_2)=s_1^*(\hat\rho_2)=t^{**}(\hat\rho_2)=\frac{1}{\mu_2}.
\EQN
Moreover, for $ \mu_1/\mu_2< \rho<1$ we have
\begin{itemize}
\item[(i).]   $t^*(\rho) <s_1^{*}(\rho) \ \text{for\ all}\ \rho\in(\mu_1/\mu_2,\hat\rho_2), \ \ \ t^*(\rho)>s_1^{*}(\rho) \ \text{for\ all}\ \rho\in(\hat\rho_2,1).$
\item[(ii).]   $t_B(\rho) <s_1^{*}(\rho) \ \text{for\ all}\ \rho\in(\mu_1/\mu_2,\hat\rho_2), \ \ \ t_B(\rho)>s_1^{*}(\rho) \ \text{for\ all}\ \rho\in(\hat\rho_2,1).$
\item[(iii).]   $t^{**}(\rho) <s_1^{*}(\rho) \ \text{for\ all}\ \rho\in(\mu_1/\mu_2,\hat\rho_2), \ \ \ t^{**}(\rho)>s_1^{*}(\rho) \ \text{for\ all}\ \rho\in(\hat\rho_2, \Lcc{\hat\rho}).$
\item[(iv).]   $t^{**}(\rho) <t^{*}(\rho) \ \text{for\ all}\ \rho\in(\mu_1/\mu_2,\hat\rho_2), \ \ \ t^{**}(\rho)>t^{*}(\rho) \ \text{for\ all}\ \rho\in(\hat\rho_2, \Lcc{\hat\rho}).$
\item[(v).]   $t^{**}(\rho) <t_B(\rho) \ \text{for\ all}\ \rho\in(\mu_1/\mu_2,\hat\rho_2), \ \ \ t^{**}(\rho)>t_B(\rho) \ \text{for\ all}\ \rho\in(\hat\rho_2, \Lcc{\hat\rho}).$
\end{itemize}

\EEL

Recall that by definition $g_L(s)=g_A(s,s)=g_B(s,s), s>0$ (cf. \eqref{eq:fL}). For the minimum of $g(t,s)$ on $\oB$ we have the following lemma.
\BEL \label{Lem:gB}
We have
\begin{itemize}
\item[(i).]   If $  0 \le \rho< \hat\rho_2$, then
\BQNY
\inf_{(t,s)\in \oB}g(t,s)=   g_L(t^{*}) =\frac{2}{1+\rho}(\mu_1+\mu_2+2/t^*),
\EQNY
where $(t^*,t^*)$ is the unique   minimizer of $g(t,s)$ on $\oB$.
\item[(ii).] If $ \rho  = \hat\rho_2$, then $t^*(\hat\rho_2)=s^*(\hat\rho_2)=1/\mu_2$, and
\BQNY
\inf_{(t,s)\in \oB}g(t,s)=   g_L({1}/{\mu_2})=g_2( {1}/{\mu_2})=4\mu_2,
\EQNY
where the  minimum  of $g(t,s)$ on $\oB$ is attained at $(1/\mu_2 ,1/\mu_2 )$,  with
$g_3(1/\mu_2 ,1/\mu_2 )=g_2(1/\mu_2)$, and $1/\mu_2 $ is the unique minimizer of $g_2(s), s\in(0,\IF).$

\item[(iii).] If  $ \hat\rho_2<\rho <1$, then  
\BQNY
 \inf_{(t,s)\in \oB}g(t,s)=   \inf_{(t,s)\in D_2}g_2(s)  =  g_2({1}/{\mu_2})=4\mu_2,
\EQNY
where
the  minimum  of $g(t,s)$ on $\oB$ is attained when $g(t,s)=g_2(s)$ on $D_2$ \Lcc{(see Figure 1)}. 
\end{itemize}
\EEL

Next we consider the minimum of $g(t,s)$ on $\oA$.  Recall $s^*(\rho)$ defined in \eqref{eq:stStar}, $s_A(\rho)$ defined in \eqref{eq:sA}, and $s^{**}(\rho)$ defined in \eqref{eq:SDoubs}. We first give the following lemma.

\BEL\label{Lem:rhoh1}
We have

(a).  Both  $s^*(\rho)$   and $s^{**}(\rho)$ are decreasing  functions on $[0,1]$.

(b).  That $\hat\rho_1$  is the unique point on $[0,1)$ such that
\BQNY
s_A(\hat\rho_1) =s^{**}(\hat\rho_1)=s^{*}(\hat\rho_1),\ \
\EQNY
and
\begin{itemize}
\item[(i).] $s_A(\rho) <s^{**}(\rho) \ \text{for\ all}\ \rho\in[0,\hat\rho_1), \ \ \ s_A(\rho)>s^{**}(\rho) \ \text{for\ all}\ \rho\in(\hat\rho_1,1)$,
\item[(ii).] $s^*(\rho) <s^{**}(\rho) \ \text{for\ all}\ \rho\in[0,\hat\rho_1), \ \ \ s^*(\rho)>s^{**}(\rho) \ \text{for\ all}\ \rho\in(\hat\rho_1,1).$

\end{itemize}
(c).  For all $\mu_1/\mu_2<\rho<1$, it holds that $s^{**}(\rho)<s_1^*(\rho)$.

\EEL


For the minimum of $g(t,s)$ on $\oA$ we have the following lemma.

\BEL \label{Lem:gA}We have
\begin{itemize}
\item[(i).] If  $0\le \rho < \hat \rho_1$, then
\BQNY
\inf_{(t,s)\in \oA}g(t,s) =   g_A(t_A,s_A)= 4(\mu_2+(1-2\rho)\mu_1),
\EQNY
where   $(t_A,s_A)\in A$ is the unique   minimizer of $g(t,s)$ on $\oA$.
\item[(ii).] If $\rho = \hat \rho_1$, then
\BQNY
\inf_{(t,s)\in \oA}g(t,s)= g_A(t_A,s_A)= 4(\mu_2+(1-2\rho)\mu_1),
\EQNY
where  $(t_A,s_A)=(t^*,s^*)\in L$ is  the unique   minimizer of $g(t,s)$ on $\oA$.
\item[(iii).]   If $ \hat \rho_1 < \rho< \hat \rho_2$, then
\BQNY
\inf_{(t,s)\in \oA}g(t,s)=  g_L(s^{*})= \frac{2}{1+\rho}(\mu_1+\mu_2+2/s^*),
\EQNY
where $(s^*,s^*)$ is the unique   minimizer of $g(t,s)$ on $\oA$.
\item[(iv).] If $ \rho  = \hat \rho_2$, then $t^*(\hat\rho_2)=s^*(\hat\rho_2)=1/\mu_2$, and
\BQNY
 \inf_{(t,s)\in \oA}g(t,s)=   g_L(s^{*})=g_2({1}/{\mu_2})=4\mu_2,
\EQNY
where the  minimum  of $g(t,s)$ on $\oA$ is attained at $(1/\mu_2 ,1/\mu_2 )$,  with
$g_3(1/\mu_2 ,1/\mu_2 )=g_2(1/\mu_2)$. 
\item[(v).] If  $\hat \rho_2<\rho <1$, then
\BQNY
\inf_{(t,s)\in \oA}g(t,s) =  g_2({1}/{\mu_2})=4\mu_2,
\EQNY
where
the  minimum  of $g(t,s)$ on $\oA$ is attained when $g(t,s)=g_2(s)$ on $D_2$ (see Figure 1). 
\end{itemize}
\EEL


Consequently, combining the results in \nelem{Lem:gB} and \nelem{Lem:gA}, we conclude that Proposition \ref{Lem:Optg} holds for $0\le \rho<1$. Thus, the proof is complete.
\QED

\section{Conclusion and discussions}\label{Sec:CD}
In the multi-dimensional risk theory, the so-called ``ruin'' can be defined
in different manner.
\kk{Motivated by diffusion approximation approach,} in this paper we modelled the risk process
via a multi-dimensional BM with drift.
\kk{We analyzed} the component-wise infinite-time ruin probability for dimension $d=2$ by solving a two-layer
optimization problem, \kk{which by the use of Theorem 1 from \cite{RolskiSPA}} led to} the logarithmic asymptotics for $P(u)$
\kk{as $u\to \infty$, given by explicit form of the adjustment coefficient $\gamma=g(\vk{t}_0)/2$ (see (\ref{eq:logri1})).
An important tool here is Lemma \ref{AL} on the quadratic programming, cited from \cite{HA2005}.}
In this way we were also able to identify the dominating points by careful analysis of different regimes for $\rho$
\kk{and specify} three regimes with different formulas for $\gamma$ (see Theorem \ref{p.log}).
An open and difficult problem is \kk{the} derivation of exact asymptotics for   $P(u)$ in \eqref{eq:PPu1},
for which the problem of finding dominating points would be the first step. A refined double-sum method as in \cite{DHJT18} might be suitable for this purpose.
A detailed analysis of the case for dimensions $d>2$ seems to be technically very complicated,
even for \kk{getting the logarithmic asymptotics.} We also note that a more \kk{natural} problem
\kk{of} 
considering $R_i(t)=\alpha_iu+\mu_it-X_i(t)$, with general $\alpha_i>0, i=1,2$,
\kk{leads to much more difficult technicalities with the analysis of $\gamma$.}

Define the ruin time of component $i$, \kk{$1\le i\le d$,} by
$T_i=\min\{t: R_i(t)<0\}$ and let $T_{(1)}\le T_{(2)}\le \ldots\le  T_{(d)}$ be \kk{the order statistics of ruin times}.
\kk{Then the component-wise infinite-time ruin probability is equivalent to $\pk{T_{(d)}<\infty}$}
\kk{while the} ruin time \kk{of at least one business line} is $T_{\min}=T_{(1)}=\min_iT_i$. Other
\kk{interesting problems}
like $\pk{T_{(j)}<\infty}$ have not yet been analysed.
\kk{For instance, it} would be interesting for $d=3$ to study the case $T_{(2)}$.
The general scheme  on  how to obtain logarithmic asymptotics for such problems
was  discussed  in \cite{RolskiSPA}.

Random vector $\vk{\bar{X}}=(\sup_{t\ge 0}(X_1(t)- p_1t),\ldots,\sup_{t\ge 0}(X_d(t)-p_dt))^{\top}$ has exponential marginals and if it is not concentrated on a subspace of dimension less than $d$, it defines a multi-variate exponential distribution. In this paper for dimension $d=2$, we derived some asymptotic properties of such distribution. Little is known about properties of this multi-variate distribution and
more studies on it would be of interest. For example a
correlation structure of $\vk{\bar{X}}$ is unknown. In particular, \kk{in the context of findings presented in this contribution,
it would be interesting to find}
the correlation
between $\sup_{t\ge0} (X_1(t)-\mu_1t)$ and $ \sup_{t\ge0} (X_2(t)-\mu_2t)$.

\section*{Appendix}

{\bf Proof of Lemma \ref{Lem:qpp}:} Referring to Lemma \ref{AL}, we have, for any fixed $t,s$, there exists a unique index set
$$
I(t,s) \subseteq \{1,2\}
$$
such that
\BQN \label{eq:gstI}
g(t,s)=(1+\mu_1 t, 1+\mu_2 s)_{I(t,s)} \ (\Sigma_{ts})_{{I(t,s)},{I(t,s)}}^{-1} \  (1+\mu_1 t, 1+\mu_2 s)_{I(t,s)}^\top,
\EQN
and
\BQN\label{eq:SigII}
(\Sigma_{ts})_{{I(t,s)},{I(t,s)}}^{-1} \  (1+\mu_1 t, 1+\mu_2 s)_{I(t,s)}^\top > \vk 0_{I(t,s)}.
\EQN
Since $I(t,s)=\{1\}, \{2\}$ or $\{1,2\}$, we have that
\begin{itemize}
\item[(S1).] On the set $E_1=\{(t,s): \rho \ t \wedge s\  s^{-1}(1+\mu_2 s)\ge (1+\mu_1 t)\}$,\ \ \  $g(t,s)=g_2(s)$
\item [(S2).] On the set $E_2=\{(t,s): \rho \ t \wedge s\  t^{-1}(1+\mu_1 t)\ge (1+\mu_2 s)\}$,  \ \ \ $g(t,s)=g_1(t)$
\item [(S3).] On the set $E_3 = (0,\IF)^2 \setminus (E_1\cup E_2)$,\ \ \ $g(t,s)=g_3(t,s)$.
\end{itemize}

Clearly, if $\rho\le 0$ then
$$
E_1=E_2=\emptyset,\ \ E_3=(0,\IF)^2.
$$
In this case,
$$g(t,s)=g_3(t,s),\ \ \ \ (t,s)\in (0,\IF)^2. $$

Next, we focus on the case where $\rho>0$. We consider the regions $\oA$ and $B$ separately.

\underline{Analysis on $\oA$.}
We have
\BQNY 
&& A_1=\oA\cap E_1 =\{s\le t \le f_1(s)\}, \ \ \ f_1(s)=\frac{\rho-1}{\mu_1}+\frac{\rho \mu_2}{\mu_1}s,\\
&& A_2=\oA\cap E_2 =\{s\le t \le f_2(s)\}, \ \ \ f_2(s)= \frac{\rho  s}{1+(\mu_2-\rho \mu_1)s}, \nonumber \\
&& A_3= \oA\cap E_3 =\{t \ge s, t>\max( f_1(s), f_2(s))\}.\nonumber
\EQNY
Next we analyse the intersection situation of the functions $f(s)=s, f_1(s),f_2(s)$ on region $\oA$.

Clearly, for any $s>0$ we have $f_2(s)<s$. Furthermore, $f_1(s)=f_2(s)$ has a unique positive solution $s_1$ given by
$$
s_1=\frac{1-\rho}{\rho (\mu_2 -\rho \mu_1)}.
$$
Finally, for $\rho \mu_2\le \mu_1$ we have
that $f_1(s)$ does not intersect with $f(s)$ on $(0,\IF)$, but for $\rho \mu_2> \mu_1$ the unique intersection point is given by  $s_1^*>s_1$ (cf. \eqref{eq:s1}).
To conclude, we have, for  $\rho\le \mu_1/\mu_2,$
$$
g(t,s)=g_3(t,s),\ \ \ (t,s)\in \oA,
$$
and for  $\rho> \mu_1/\mu_2,$
$$
g(t,s)=\left\{\begin{array}{cc}
g_3(t,s), & \hbox{if }  (t,s)\in  \oA\cap \{t\ge \max(s,f_1(s)), t> f_1(s)\}\\
 g_2(s), & \hbox{if }   (t,s)\in  \oA\cap\{s\le t\le f_1(s)\}.
\end{array}\right.
$$
Additionally, we have from Lemma \ref{AL} $g_3(f_1(s), s)= g_2(s)$ for all $s\ge s_1^*$.

\underline{Analysis on $B$.} The two scenarios $\rho\le \mu_1/\mu_2$ and $\rho> \mu_1/\mu_2$ will be considered separately.
 For $\rho\le \mu_1/\mu_2,$  we have
\BQNY
&& B_1=B\cap E_1 =\{ t< s \le h_1(t)\}, \ \ \ h_1(t)= \frac{\rho  t}{1+(\mu_1-\rho \mu_2)t}, \\
&& B_2=B\cap E_2 =\{t <s \le h_2(t)\}, \ \ \ h_2(t)= \frac{\rho-1}{\mu_2}+\frac{\rho \mu_1}{\mu_2} t, \\
&& B_3= B\cap E_3 =\{s >\max(t,  h_1(t), h_2(t))\}.
\EQNY
It is easy to check that
$$
t>h_1(t),\ \ \ t>h_2(t),\ \ \ \forall t>0,
$$
and thus
$$
g(t,s)=g_3(t,s),\ \ \ (t,s)\in B.
$$
For $\rho> \mu_1/\mu_2,$ we have
\BQNY
&& B_1=B\cap E_1 =\{ w_1(s)\le t< s  \}, \ \ \ w_1(s)= \frac{s}{\rho+( \rho \mu_2-\mu_1)s},\\
&& B_2=B\cap E_2 =\{w_2(s)\le t <s \}, \ \ \ w_2(s)= \frac{ 1-\rho}{\mu_1\rho }+\frac{  \mu_2}{\mu_1\rho} s, \nonumber\\
&& B_3= B\cap E_3 =\{ t<\min(s, w_1(s), w_2(s))\}.\nonumber
\EQNY

Next we analyze the intersection situation of the functions $w(s)=s, w_1(s),w_2(s)$ on region $B$.

Clearly, for any $s>0,$ $w_2(s)>s$. $w_1(s)$ and $w_2(s)$ do not intersect on $(0,\IF)$. $w(s)$ and $w_1(s)$ has a unique intersection point $s_1^*$ (cf. \eqref{eq:s1}).

To conclude, we have, for $\rho\le \mu_1/\mu_2,$
$$
g(t,s)=g_3(t,s),\ \ \ (t,s)\in B,
$$
and for  $\rho> \mu_1/\mu_2,$
$$
g(t,s)=\left\{\begin{array}{cc}
g_3(t,s), & \hbox{if }  (t,s)\in  B\cap \{t< \min(s,w_1(s))\}\\
 g_2(s), & \hbox{if }   (t,s)\in  B\cap\{ w_1(s)\le t<s  \}.
\end{array}\right.
$$
Additionally, it follows from Lemma \ref{AL} that $g_3(w_1(s), s)= g_2(s)$ for all $s\ge s_1^*$.

Consequently, the claim follows by a combination of the above results. This completes the proof.
\QED

 {\bf Proof of Lemma \ref{Lem:tt}.}  (a). The claim for $t^*(\rho)$ follows by noting its following representation:
\BQNY
t^*(\rho)=s^*(\rho)=\sqrt {\frac{2(1-\rho)}{\mu_1^2+\mu_2^2-2\mu_1\mu_2+2\mu_1\mu_2-2\rho\mu_1\mu_2}}=\sqrt { \frac{2}{\frac{(\mu_1-\mu_2)^2}{1-\rho}+2\mu_1\mu_2} }.
\EQNY
The claims for $t_B(\rho)$ and $s_1^*(\rho)$ follow  directly from their definition.

(b).
First note that
\BQNY
t^{**}(0)=t^{**}( \hat\rho_2)=\frac{1}{\mu_2}.
\EQNY
Next it is calculated that
\BQNY
 \frac{\partial t^{**}(\rho)}{\partial \rho} =\frac{-2\mu_2 \rho^2+4\mu_2\rho -\mu_1-\mu_2}{(\mu_2+\rho \mu_1-2\mu_2\rho^2)^2}.
\EQNY
Thus, the claim of (b) follows by analysing the sign of $\frac{\partial t^{**}(\rho)}{\partial \rho} $ over $(0,\hat \rho)$.

(c).
For any $0\le \rho\le \mu_1/\mu_2$ we have $\abs{\mu_1-2\rho\mu_2}\le \mu_1$ and thus
\BQNY
t_B(\rho)\ge \frac{1}{u_1}\ge\frac{1}{u_2}\ge\frac{1-\rho}{u_2(1-\rho^2)}\ge \frac{1-\rho}{\rho(\mu_1-\rho \mu_2)+\mu_2(1-\rho^2)}=t^{**}(\rho).
\EQNY
Further, since
\BQNY
\mu_1^2+\mu_2^2-2\rho\mu_1\mu_2=\mu_1(\mu_1-\rho\mu_2)+\mu_2(\mu_2-\rho\mu_1)\le \mu_2(\mu_1-\rho\mu_2)+\mu_2(\mu_2-\rho\mu_1) \le 2\mu_2^2(1-\rho),
\EQNY
it follows that
\BQNY
t^*(\rho)\ge \frac{1}{\mu_2}\ge t^{**}(\rho).
\EQNY
(d).  It is easy to check that \eqref{eq:tttt} holds.   For (i) we have
\BQNY
t^*(\rho) -s_1^{*}(\rho) =(1-\rho)\LT( \frac{1}{ f_1(\rho)} -\frac{1}{f_2(\rho)}\RT),
\EQNY
where
\BQNY
f_1(\rho)&=&\sqrt{\frac{(1-\rho) (\mu_1^2+\mu_2^2-2\rho\mu_1\mu_2)}{2}}=\sqrt{\mu_1\mu_2\rho^2-\frac{(\mu_1+\mu_2)^2}{2}\rho +\frac{\mu_1^2+\mu_2^2}{2}}\\
f_2(\rho)&=&\rho\mu_2-\mu_1.
\EQNY
Analysing the properties of the above two functions, we have $f_1(\rho)$ is strictly decreasing on $[0,1]$ with
$$
f_1(0)=\sqrt{\frac{\mu_1^2+\mu_2^2}{2}}>-\mu_1=f_2(0),\ \ \ f_1(1)=0\le \mu_2-\mu_1=f_2(1),
$$
and thus there is a unique intersection point of the two curves $t^*(\rho)$ and $s_1^{*}(\rho) $ which is  $\rho=\hat\rho_2$. Therefore, the claim of (i) follows.
Similarly, the claim of (ii) follows since
\BQNY
t_B(\rho)-s_1^*(\rho)=\frac{-(\mu_1+\mu_2)\rho +2\mu_2\rho^2}{(\rho\mu_2-\mu_1)(2\rho\mu_2-\mu_1)}.
\EQNY
 Finally, the claims of (iii), (iv) and (v) follow easily from (a), (b) and \eqref{eq:tttt}. This completes the proof. \QED

{\bf Proof of Lemma \ref{Lem:gB}.} Consider first the case where $0\le \rho \le \mu_1/\mu_2$. Recall \eqref{eq:ss}. We check if any of $s_i, i=1,2,$ is greater than $t$. Clearly, $s_1\le t$. Next, we check 
 whether $s_2>t$. 
It is easy to check that
\BQNY
s_2>t \ \ \Leftrightarrow \ \  t< t^{**},
\EQNY
where \Lcc{(recall \eqref{eq:TDoubs})}
\BQNY
t^{**}=t^{**}(\rho)=\frac{1-\rho}{\rho(\mu_1-\mu_2\rho)+\mu_2(1-\rho^2)}>0.
\EQNY
Then
$$
\inf_{(t,s)\in \oB} g_3(t,s) =\min\LT(\inf_{0<t<t^{**}} g_B(t,s_2(t)),  \inf_{t\ge t^{**}} g_B(t,t)\RT).
$$

Consequently, it follows from (c) of Lemma \ref{Lem:tt} the claim of (i) holds for $0\le \rho\le {\mu_1}/{\mu_2}$.

Next, we consider $\mu_1/\mu_2<\rho<1$.
Recall the function $w_1(s)$ defined in \eqref{eq:fw}. Denote the inverse function of $w_1(s)$ by
\BQNY
\hat w_1(t)=\frac{\rho t}{1-(\rho\mu_2-\mu_1) t}, \ \ t\ge s_1^*.
\EQNY
We have from Lemma \ref{Lem:qpp} that
\BQNY
 g_B(t,\hat w_1(t)) = g_2(t),\ \ \ t\ge s_1^*.
\EQNY
Further note that ${1}/{\mu_2}$ is the unique minimizer of $g_2(s), s> 0$. For $\mu_1/\mu_2<\rho< \hat\rho_2,$ we have from   (d) in Lemma \ref{Lem:tt}  that 
\BQNY
\inf_{s_1^*\le s}g_2(s) = g_2(s_1^*) =g_L(s_1^*) >g_L(t^*),
\EQNY
and further
\BQNY
\inf_{(t,s)\in \oB}g(t,s)&=&\min( \inf_{0<t<t^{**}}g_B(t,s_2(t)),  \inf_{t^{**}\le t<s_1^*}g_B(t,t),  \inf_{s_1^{*}\le t}g_B(t,\hat w_1(t)),  \inf_{s_1^*\le s}g_2(s)  )\\
&=&g_B(t^*,t^*)=g_L(t^*),
\EQNY
where $(t^*,t^*)$ is the unique minimizer of $g(t,s)$ on $\oB$. Therefore, the claim for $\mu_1/\mu_2<\rho< \hat\rho_2$ is established.

 For $\rho=\hat\rho_2$, because of \eqref{eq:tttt} we have
\BQNY
\inf_{(t,s)\in \oB}g(t,s)&=&\min( \inf_{0<t<1/\mu_2}g_B(t,s_2(t)), \inf_{ 1/\mu_2\le t}g_B(t,\hat w_1(t)),  \inf_{1/\mu_2\le s}g_2(s)  )\\
&=& g_B(1/\mu_2,1/\mu_2) =g_L(1/\mu_2)=g_2(1/\mu_2),
\EQNY
and the unique minimum  of $g(t,s)$ on $\oB$ is attained at $(1/\mu_2 ,1/\mu_2 )$. 
Moreover, for all $\hat\rho_2<\rho<  1$ we have
\BQNY
s_2(t_B)=\hat w_1(t_B)=\frac{1}{\mu_2}>s_1^*.
\EQNY
Thus,
\BQNY
\inf_{(t,s)\in \oB}g(t,s)&=&\min( \inf_{0<t<t_B}g_B(t,s_2(t)),  \inf_{t_B\le t}g_B(t,\hat w_1(t)),  \inf_{s_1^*\le s}g_2(s)  )\\
&=& g_B(t_B,  {1}/{\mu_2})= g_2( {1}/{\mu_2}),
\EQNY
and the unique minimum  of $g(t,s)$ on $\oB$ is attained when $g(t,s)=g_2(s)$ on $D_2$.
This completes the proof. \QED

{\bf Proof of Lemma \ref{Lem:rhoh1}.} 
(a). The claim for $s^*(\rho)$ has been shown in the proof of (a) in Lemma \ref{Lem:tt}.
Next, we show the claim for $s^{**}(\rho)$, for which it is sufficient to show that $\frac{\partial s^{**}(\rho)}{\partial \rho}<0$ for all $\rho\in[0,1]$. In fact, we have
\BQNY
\frac{\partial s^{**}(\rho)}{\partial \rho} =\frac{-2\mu_1 \rho^2+4\mu_1\rho -\mu_1-\mu_2}{(\mu_1+\rho \mu_2-2\mu_1\rho^2)^2}<0.
\EQNY
(b).  In order to prove (i), the following two scenarios will be discussed separately:
\BQNY
(S1).\ \ \ \mu_2<2\mu_1; \ \ \ \ \ \ (S2).\ \ \ \mu_2\ge2\mu_1.
\EQNY
First consider (S1).
If $0\le\rho<\frac{\mu_2}{2\mu_1}$, then
\BQNY
s_A(\rho)-s^{**}(\rho)&=&\frac{(\mu_1+\rho\mu_2-2\mu_1\rho^2)-(1-\rho)(\mu_2-2\rho\mu_1)}{(\mu_2-2\rho\mu_1)(\mu_1+\rho\mu_2-2\mu_1\rho^2)}\\
&=& \frac{f(\rho)}{(\mu_2-2\rho\mu_1)(\mu_1+\rho\mu_2-2\mu_1\rho^2)},
\EQNY
where
\BQNY
f(\rho)=- 4\mu_1\rho^2 + 2(\mu_2+\mu_1)\rho- \mu_2+\mu_1.
\EQNY
Analysing the function $f$, we conclude that
\BQNY
f(\rho)<0,\ \ \ \text{for}\ \ \rho\in[0,\hat\rho_1),\ \ \ \ f(\rho)  >  0,\ \ \ \text{for}\ \ \rho\in(\hat\rho_1, \frac{\mu_2}{2\mu_1}).
\EQNY
Further, for $\frac{\mu_2}{2\mu_1}\le \rho<1 $ we have
\BQNY
s_A(\rho)-s^{**}(\rho) = \frac{\mu_1+\mu_2-2\mu_1\rho}{( 2\rho\mu_1-\mu_2)(\mu_1+\rho\mu_2-2\mu_1\rho^2)} >0.
\EQNY
Thus, the claim in (i) is established for (S1).  Similarly, the claim in (i) is valid for (S2) .
 Next, note that
\BQNY
s^*( \rho ) -s^{**}( \rho )=(1-\rho)\LT(\frac{1}{f_1(\rho)}-\frac{1}{f_2(\rho)}\RT)
\EQNY
with
\BQNY
f_1(\rho)&=&\sqrt{\frac{(1-\rho) (\mu_1^2+\mu_2^2-2\rho\mu_1\mu_2)}{2}}=\sqrt{\mu_1\mu_2\rho^2-\frac{(\mu_1+\mu_2)^2}{2}\rho +\frac{\mu_1^2+\mu_2^2}{2}}\\
f_2(\rho)&=&\mu_1+\rho \mu_2-2\mu_1\rho^2.
\EQNY
Analysing the properties of the above two functions, we have $f_1(\rho)$ is strictly decreasing on $[0,1]$ with
$$
f_1(0)=\sqrt{\frac{\mu_1^2+\mu_2^2}{2}}\ge \mu_1=f_2(0),\ \ \ f_1(1)=0\le \mu_2-\mu_1=f_2(1),
$$
and thus there is a unique intersection point $\rho\in (0,1)$ of $s^*( \rho )$ and $s^{**}( \rho )$. It seems not clear at the moment whether this unique point is $\hat\rho_1$ or not, since we have to solve a polynomial equation of order 4.  Instead, it is sufficient to show that
\BQN\label{eq:sAsS}
s_A(\hat\rho_1)=s^*(\hat\rho_1).
\EQN
In fact, basic calculations show that the above is equivalent to
\BQNY
(2\mu_1 \hat\rho_1 - (u_1+\mu_2)) f(\hat\rho_1)=0
\EQNY
which is valid due to the fact that $f(\hat\rho_1)=0$.
Finally, the claim in (c) follows since
$$
\rho \mu_2-\mu_1<\mu_1+\rho\mu_2-2\rho^2\mu_1.
$$
This completes the proof. \QED

{\bf Proof of Lemma \ref{Lem:gA}.}  Two cases
$\hat\rho_1\le \mu_1/\mu_2$ and $\hat\rho_1>\mu_1/\mu_2$
should be distinguished.
Since the proofs for these two cases are similar, we give below only the proof for the more complicated case $\hat\rho_1\le \mu_1/\mu_2$.

Note that, for $0\le \rho\le \mu_1/\mu_2$, as in \eqref{eq:g3A},
$$
\inf_{(t,s)\in \oA} g(t,s) =\inf_{(t,s)\in \oA} g_3(t,s)=\min\LT(\inf_{0<s<s^{**}} f_A(s),  \inf_{s\ge s^{**}} g_L(s)\RT),
$$
and thus the claim for $0\le \rho\le \mu_1/\mu_2$  follows directly from (i)-(ii) of (b) in Lemma \ref{Lem:rhoh1}.
Next, we consider the case $\mu_1/\mu_2<\rho<\hat\rho_2$ (note here $\hat\rho_1<\mu_1/\mu_2<\rho$). We have by (i) of (d) in \nelem{Lem:tt} and (i)-(ii) of (b) in Lemma \ref{Lem:rhoh1} that
 \BQNY 
 s^{**}(\rho)<s^*(\rho)=t^*(\rho) <s_1^{*}(\rho), \ \ \ s_1^{*}(\rho)>\frac{1}{\mu_2},\ \ \ s_A(\rho)>s^{**}(\rho).
 \EQNY
Thus, it follows from Lemma \ref{Lem:qpp} that
\BQNY
\inf_{(t,s)\in \oA} g(t,s)&=&\min\LT(\inf_{0<s<s^{**}} g_A(t_2(s),s),  \inf_{ s^{**}\le s \le s_1^*} g_A(s,s), \inf_{s_1^* <s} g_A(f_1(s),s), \inf_{{s_1^* <s} } g_2(s)\RT) \\
&=& g_A(t^*,s^*)=g_L(s^*),
\EQNY
and $(t^*,s^*)\in L$ is the unique minimizer of $g(t,s)$ on $\oA$.
Here we used the fact that 
\BQNY
\inf_{s_1^* <s} g_A(f_1(s),s) =\inf_{{s_1^* <s} } g_2(s)=g_A(f_1(s_1^*),s_1^*) =g_2(s_1^*)>g_L(s^*).
\EQNY
Next, if $\rho=\hat\rho_2$, then
\BQNY
s_1^*(\hat\rho_2)=s^*(\hat\rho_2)=\frac{1}{\mu_2},
\EQNY
and thus
\BQNY
\inf_{(t,s)\in \oA} g(t,s) &=&\min\LT(\inf_{0<s<s^{**}} g_A(t_2(s),s),  \inf_{ s^{**}\le s \le  {1}/{\mu_2}} g_A(s,s), \inf_{ {1}/{\mu_2} <s} g_A(f_1(s),s), \inf_{{{1}/{\mu_2} <s} } g_2(s)\RT)\\
&= & g_A( 1/\mu_2, 1/\mu_2 )=g_L( 1/\mu_2)=g_2( 1/\mu_2).
\EQNY
Furthermore,  the unique minimum  of $g(t,s)$ on $A$ is attained at $(1/\mu_2 ,1/\mu_2 )$,  with
$g_3(1/\mu_2 ,1/\mu_2 )=g_2(1/\mu_2)$. 

Finally, for $\hat\rho_2<\rho<1$, we have
 \BQNY 
 s^{**}(\rho)<s_1^*(\rho)  <s^{*}(\rho)<\frac{1}{\mu_2},\ \ \ s_A(\rho)>s^{**}(\rho),
 \EQNY
and thus
\BQNY
\inf_{(t,s)\in \oA} g(t,s) &=&\min\LT(\inf_{0<s<s^{**}} g_A(t_2(s),s),  \inf_{ s^{**}\le s \le s_1^*} g_A(s,s), \inf_{s_1^* <s} g_A(f_1(s),s), \inf_{{s_1^*<s} } g_2(s)\RT)\\
&=& g_2( {1}/{\mu_2}),
\EQNY
where
the unique minimum  of $g(t,s)$ on $\oA$ is attained when $g_3(t,s)=g_2(s)$ on $D_2$. 
This completes the proof. \QED


\bigskip
{\bf Acknowledgement}: We are thankful to the referees for their carefully reading and constructive suggestions which significantly improved the manuscript.
TR \& KD acknowledge partial support by NCN Grant No 2015/17/B/ST1/01102 (2016-2019).

\bibliographystyle{plain}

 \bibliography{vectProcEKEE}
\end{document}